\documentclass[11pt,here]{article}
\usepackage[colorlinks=true,
linkcolor=webgreen,
filecolor=webbrown,
citecolor=webgreen]{hyperref}

\usepackage{amssymb,psfig,epsfig,latexsym,graphicx}

\definecolor{webgreen}{rgb}{0,.5,0}
\definecolor{webbrown}{rgb}{.6,0,0}

\newtheorem{lemma}{Lemma}
\newtheorem{theorem}{Theorem}
\newtheorem{coro}{Corollary}

\newcommand{\eqn}[1]{(\ref{#1})}
\newcommand{\bsq}{{\vrule height .9ex width .8ex depth -.1ex }}

\newcommand{\FF}{{\mathbb F}}
\newcommand{\RR}{{\mathbb R}}

\newcommand{\sD}{{\mathcal D}}

\newcommand{\sC}{{\mathcal C}}
\newcommand{\eeq}{\end{equation}}
\newcommand{\beql}[1]{\begin{equation}\label{#1}}

\makeatletter
\def\@sect#1#2#3#4#5#6[#7]#8{\ifnum #2>\c@secnumdepth
     \def\@svsec{}\else
     \refstepcounter{#1}\edef\@svsec{\csname the#1\endcsname.\hskip .75em }\fi
     \@tempskipa #5\relax
      \ifdim \@tempskipa>\z@
        \begingroup #6\relax
          \@hangfrom{\hskip #3\relax\@svsec}{\interlinepenalty \@M #8\par}%
        \endgroup
       \csname #1mark\endcsname{#7}\addcontentsline
         {toc}{#1}{\ifnum #2>\c@secnumdepth \else
                      \protect\numberline{\csname the#1\endcsname}\fi
                    #7}\else
        \def\@svsechd{#6\hskip #3\@svsec #8\csname #1mark\endcsname
                      {#7}\addcontentsline
                           {toc}{#1}{\ifnum #2>\c@secnumdepth \else
                             \protect\numberline{\csname the#1\endcsname}\fi
                       #7}}\fi
     \@xsect{#5}}
\def\@begintheorem#1#2{\it \trivlist \item[\hskip \labelsep{\bf #1\ #2.}]}

\def\section{\@startsection {section}{1}{\z@}{-3.5ex plus -1ex minus
 -.2ex}{2.3ex plus .2ex}{\normalsize\bf}}
\makeatother
\makeatletter
\def\subsection{\@startsection {subsection}{1}{\z@}{-3.5ex plus -1ex minus
 -.2ex}{2.3ex plus .2ex}{\normalsize\bf}}

\makeatother

\begin{document}
\begin{center}
{\large\bf On Asymmetric Coverings and Covering Numbers} \\
\vspace*{+.5in}
{\em David Applegate, E. M. Rains}\footnote{Present address: Center for Communications Research, Princeton, NJ 08540} and {\em N. J. A. Sloane} \smallskip \\
Information Sciences Research Center \\
AT\&T Shannon Lab \\
Florham Park, NJ 07932--0971 \medskip \\
Email addresses: \href{mailto:david@research.att.com}{{\tt david@research.att.com}}, \href{mailto:rains@idaccr.org}{{\tt rains@idaccr.org}}, \href{mailto:njas@research.att.com}{{\tt njas@research.att.com}}
\bigskip \\

{\bf Abstract}
\end{center}

An asymmetric covering $\sD (n,R)$ is a collection of
special subsets $S$ of an $n$-set such that every subset $T$
of the $n$-set is contained in at least one special $S$ with $|S| - |T| \le R$.
In this paper we compute the smallest size of any $\sD (n,1)$ for $n \le 8$.
We also investigate ``continuous'' and ``banded'' versions of the problem.
The latter involves the classical covering numbers $C(n,k,k-1)$,
and we determine the following new values:
$C(10,5,4) = 51$,
$C(11,7,6,) =84$,
$C(12,8,7) = 126$, $C(13,9,8)= 185$ and
$C(14,10,9) = 259$.
We also find the number of nonisomorphic minimal
covering designs in several cases.

%
% Section 1
%

\section{Introduction}
Let $D(n,R)$ denote the smallest size
of any asymmetric covering\footnote{These are called
{\em directed coverings} in \cite{CEK}. However, that term has already been used
in the literature with a different meaning (cf. \cite{MM92}).}
$\sD (n,R)$.
Prompted by applications to the manufacture of semiconductor wafers,
Cooper, Ellis and Kahng \cite{CEK} have investigated 
the asymptotic behavior of $D(n,R)$ for fixed $R$ as $n \to \infty$.
In Section 2 of the present paper we show that the values of $D(n,1)$ for $n \le 8$ are as shown in Table \ref{T1}.\footnote{This is sequence A66000 in \cite{OEIS}.
If any further values are computed they will be recorded there.}

Our method of attack is to formulate $D(n,1)$ as the solution to a $\{0,1\}$-integer programming problem.
If instead we allow the variables to take any real values in the range $[0,1]$ then the linear program can be solved exactly (Section 3).
Of course this provides a lower bound $E(n)$ to $D(n,1)$.

An upper bound can be obtained by restricting to asymmetric coverings with a certain banded structure defined in Section 4.
Corollary \ref{Cor2} shows that the solution $C(n)$ to the banded version of the problem is given by
\beql{EqC1}
C(n) = \sum_{i=0}^{[n/2]} C(n+1, n+1-2i, n-2i ) ~,
\eeq
where as usual $C(v,k,t)$ denotes the smallest 
size of any {\em covering design} $\sC (v,k,t)$;
that is, any collection of special $k$-subsets $S$ of a $v$-set such that any $t$-subset $T$ is contained in at least one $S$.

Although there have been a large number of papers written 
about covering designs (for recent work see \cite{CKRM}, 
\cite{Gir90}, \cite{LJCR},
\cite{MM92}, \cite{NO99}, \cite{NO99a}, \cite{CRC}),
not many exact values are known.
In Section 5 we determine several new values of $C(n,k,k-1)$.
Call a covering design $\sC (n,k,t)$ {\em optimal} if it contains the smallest number $C(n,k,t)$ of subsets, and {\em minimal} if it is
no longer a covering if any subset is omitted.
In Section 5 we also determine the number of nonisomorphic optimal covering designs $\sC (n,k,k-1)$ for $n$ up through 10,
as well as
the number of minimal covering designs for $n$ up through 7.
See Tables \ref{TC1}, \ref{TC2} in Section 5.

Using these results we obtain the values of $C(n)$ for $n \le 11$ shown in Table \ref{T1}.

Some of the literature on covering designs works with the complements of the special sets, in which case this is called the
{\em Tur\'{a}n design problem} (cf. \cite{CKW}), and the {\em Tur\'{a}n number}
$T(v,k,t)$ is equal to $C(v, v-t, v-k)$.
Of course our results also provide new values for certain Tur\'{a}n numbers.

\begin{table}[htb]
\caption{Size $D(n,1)$ of smallest asymmetric 1-covering of an $n$-set, together with values of the continuous and banded solutions $E(n)$ and $C(n)$.}
$$
\begin{array}{|c|ccc|} \hline
n & E(n) & D(n,1) & C(n) \\ \hline
1 & 1 & 1 & 1 \\
2 & 2 & 2 & 2 \\
3 & 3 & 3 & 3 \\
4 & 5 & 6 & 6 \\ [+.05in]
5 & 8 \frac{1}{2} & 10 & 10 \\ [+.05in]
6 & 14 \frac{5}{6} & 18 & 18 \\ [+.05in]
7 & 26 \frac{3}{8} & 31 & 31 \\ [+.05in]
8 & 47 \frac{23}{40} & 58 & 60 \\ [+.05in]
9 & 86 \frac{553}{720} & ? & 106 \\ [+.05in]
10 & 159 \frac{353}{560} & ? & 196 \\ [+.05in]
11 & 295 \frac{3337}{4480} & ? & 352 \\ \hline
\end{array}
$$
\label{T1}
\end{table}

\subsection*{Notation}
Let $\FF_2^n$ denote the set of binary vectors of length $n$.
We represent subsets of an $n$-set by their indicator vectors in $\FF_2^n$, and
then an asymmetric covering $\sD (n,R)$ can be thought of as a binary code called
an ``asymmetric covering code''.
As usual weight (denoted wt) and distance (dist)
refer to Hamming weight and Hamming distance.
The {\em co-weight} of $u \in \FF_2^n$ is $n-wt (u)$.
Two codes, coverings or designs are isomorphic if they differ just by a permutation of the coordinates.

\section{Values of $D(n,1)$ for $n \le 8$.}
Let $x_u$, $u \in \FF_2^n$, be real $\{0,1\}$-valued variables.
Then $D(n,1)$ is equal to the minimal value of
\beql{EqE1}
\sum_{u \in \FF_2^n} x_u
\eeq
subject to the constraints
\beql{EqE2}
\sum_{u \in \FF_2^n \atop u \supset v}
x_u ~ + ~ x_v ~ \ge ~ 1 , ~
\mbox{for all}~ v \in \FF_2^n ~,
\eeq
where $u \supset v$ indicates that $u$ covers $v$ and $dist (u,v) =1$.
If $u= 11 \ldots 1$ then necessarily $x_u =1$.
The corresponding asymmetric covering code consists of the vectors
$u \in \FF_2^n$ for which $x_u =1$.
For example $D(3,1) =3$, and the code (which is unique up to permutation of the coordinates) is
\beql{Eq2a}
\{111,110, 001\} ~.
\eeq
Every binary vector of length 3 is either in this code or is contained in a 
codeword at distance 1 below it.

We call the above minimization problem the exact integer programming (or IP) problem.
If we relax the constraints and allow the $x_u$ to take any real values in the range $[0,1]$ we get a continuous linear programming (or LP) problem, whose solution we denote by $E(n)$.

\begin{theorem}\label{thD}
The values of $D(n,1)$ for $n \le 8$ are as shown in Table \ref{T1}.
\end{theorem}

\paragraph{Proof.}
We attacked the IP problem using CPLEX \cite{CPLEX} with AMPL \cite{FGK93} as a convenient interface.
CPLEX uses a branch and bound strategy for such problems.
We regard solutions obtained in this way as perfectly rigorous, since the computations could in principle be replaced by extremely tedious hand calculations.

For $n \le 7$ CPLEX was able to find solutions directly, without any additional assumptions being added.
Explicit solutions are described in Section 4.

For $n=8$ we must show that $D(8,1) = 58$.
A solution of size 58 found by CPLEX is given in Table \ref{T2}.
(Each vector is represented by two hexadecimal characters.
This covering has no apparent structure--in particular it has
trivial automorphism group)
To show that 57 is impossible we argue as follows.
If we add the extra assumption that there are at most eight codewords of weight 5 to the continuous LP problem, the solution is at least 63.
Therefore there must be at least nine codewords of weight 5. From
the tables of constant weight codes \cite{Andw} it follows that there must be two codewords of weight 5 and distance
exactly two apart.
Without loss of generality we can assume that $u_1 = 11111000$ and $u_2 = 11110100$ are in the code.

Suppose the code contains a vector $u_3$ of weight 3 with $dist (u_1, u_3) =6$ and $dist (u_2, u_3 ) =4$, say $u_3 =00010101$.
Then CPLEX finds that the minimal solution to the IP problem is 58.
On the other hand if no such vector $u_3$ is present (this rules out 16 vectors of weight 3) no feasible solution to the IP problem of size $\le 57$ exists.
Hence $D(8,1) = 58$.~~~$\bsq$

The total computing time for these calculations was less than 48 hours.

\begin{table}[htb]
\caption{A minimal asymmetric covering $\sD (8,1)$ of length 8 containing 58 sets (represented in hexadecimal).}
\begin{center}
\begin{tabular}{c@{~}c@{~}c@{~}c@{~}c@{~}c@{~}c@{~}c@{~}c@{~}c@{~}c@{~}c@{~}c@{~}c@{~}c@{~}c@{~}c@{~}c@{~}c@{~}c}
01 & 07 & 0A & 11 & 1E & 28 & 2D & 33 & 34 & 37 & 3B & 4B & 4C & 52 & 55 & 57 & 5D & 61 & 66 & 6E \\
6F & 73 & 75 & 78 & 7E & 7F & 84 & 89 & 8F & 96 & 98 & 99 & 9F & A2 & A5 & AA & B3 & BB & BC & BD \\
C0 & C3 & CC & D5 & DA & DB & DD & E6 & E7 & E9 & EE & EF & F0 & F6 & F7 & F9 & FE & FF
\end{tabular}
\end{center}
\label{T2}
\end{table}

\section{Solution to the continuous linear programming problem}
\begin{theorem}\label{thE}
The optimal solution to the continuous LP problem of choosing $0 \le x_u \le 1$ for $u \in \FF_2^n$ so as to minimize \eqn{EqE1} subject to \eqn{EqE2} is given by
\beql{EqE3}
E(n) = (-1)^n n! \{ R_n (2) - R_n (1) R_{n-1} (1) \}
\eeq
where
\beql{EqE4}
R_n (x) = \sum_{k=0}^n \frac{(-x)^k}{k!}
\eeq
is the degree $n$ partial sum of $e^{-x}$.
\end{theorem}

\paragraph{Proof.}
We may assume that $x_u$ depends only on the weight of $u$.
(For let $y_u$ denote the average value of $x_v$ over all $v$ with $wt (v) = wt (u)$.
Then by averaging \eqn{EqE2} we see that the $y_u$ satisfy the same constraints as the $x_u$, and
$$\sum_{u \in \FF_2^n} y_u = \sum_{u \in \FF_2^n} x_u ~.
$$
So a symmetrized solution is just as good as a general solution.)

The ``weight enumerator'' of a symmetrized solution $x_u$ is defined by
$A_w = \sum_{wt(u) =w} x_u$ for $w=0,\ldots, n$.
The quotes are needed because the $x_u$ are in general not integers.
Let
$A(z) = \sum_{w=0}^n A_w z^w$.

The covering condition \eqn{EqE2} reads
$$wA_w + A_{w-1} \ge \left( \begin{array}{c}
n \\ w-1 \end{array}\right) , \quad w=1,\ldots, n ~,
$$
or in other words
\beql{EqE5}
A(z) + A' (z) \ge (z+1)^n ~.
\eeq
We wish to choose $A_0 , \ldots, A_n \ge 0$ so as to minimize $A(1)$ subject to
\eqn{EqE5}.
The dual problem (compare \cite[Chapt. 17]{MS77}) is to choose $B_0, \ldots, B_n \ge 0$ so as to maximize $B(1)$ subject to
\beql{EqE6}
B(z) + B' (z) \le (z+1)^n
\eeq
where $B(z) = \sum_{w=0}^n B_w z^w$.
We claim that
\beql{EqE7}
B(z) = (-1)^n n! \{ R_n (z+1) - R_{n-1} (1) R_n (z) \}
\eeq
is a feasible solution to the dual problem.
In fact it is straightforward to verify that $B_w \ge 0$ for all $w$ and
$$B(z) + B' (z) = (z+1)^n - R_{n-1} (1) z^n ~.$$
Since $R_{n-1} (1) \ge 0$,
\eqn{EqE6} holds.

Therefore
$$B(1) = (-1)^n n! \{R_n (2) - R_n (1) R_{n-1} (1) \}
$$
is an upper bound to the optimal solution to the primal problem.

On the other hand
$$A(z) = (-1)^n n! \{R_n (z+1) - R_n (1) R_{n-1} (z) \}$$
satisfies
$$A(z) + A' (z) = (z+1)^n + R_n (1) nz^{n-1}$$
and is easily checked to be a feasible solution to the primal problem.
Since $A(1) = B(1)$, this must be the optimal solution to both problems.~~~$\bsq$

\begin{coro}\label{CorA}
As $n \to \infty$,
\beql{EqE99}
E(n) \sim 2^{n+1} \left( \frac{1}{n} - \frac{3}{n^2} + O \left( \frac{1}{n^3} \right) \right) ~.
\eeq
\end{coro}
We omit the routine derivation of this from \eqn{EqE3}.

The first few values of $E(n)$ are shown in Table \ref{T1}.

\section{Banded solutions}
Let $\sD$ be an asymmetric covering $\sD (n,1)$.
We call $\sD$ {\em banded} if every vector $v \in \FF_2^n$ with odd
co-weight is covered by a vector $u \in \sD$ of weight one higher.

For example \eqn{Eq2a} is banded, since the vector 000 is covered by 001
and the vectors 011, 101, 110 are all covered by 111.

\begin{theorem}\label{thBC}
If a code $\sC \subseteq \FF_2^{n+1}$ is a union of covering designs,
\beql{EqE8}
\sC = \bigcup_{i=0}^{[n/2]} \sC (n+1, n+1-2i, n-2i ) ~,
\eeq
then deleting\footnote{Or puncturing, cf. \cite{MS77}, p.28.} any one
coordinate from all the vectors of $\sC$ yields a banded asymmetric covering $\sD (n,1)$.
Conversely, let $\sD$ be a banded asymmetric covering $\sD (n,1)$.
If we append a $0$ or $1$ to every vector of $\sD$ in such a way that all co-weights become even, the result is a union of covering designs of the form \eqn{EqE8}.
\end{theorem}

\paragraph{Proof.}
Suppose $\sC$ has the structure shown in \eqn{EqE8} and let $\sD$ be obtained by deleting one coordinate, which for concreteness we suppose is the last
coordinate.
We must show that $\sD$ is a banded asymmetric covering.
Let $v \in \FF_2^n$ have weight $w$.
If the co-weight $n-w$ is even, say $2i$, then $v^\ast = v0$ must be covered by some vector $u^\ast = u \delta$, $\delta =0$ or 1, in the covering design
$\sC (n+1, n+1-2i, n-2i )$, and then $u \in \sD$ covers $v$.
On the other hand if $n-w$ is odd, say $2i -1$, then $v^\ast = v1$ must be covered by some $u^\ast = u1 \in \sC (n+1, n+1-2i, n-2i )$, and again $u$ covers $v$.
The converse is established by similar arguments.

Since the covering number $C(n,k,v)$ is by definition the size of the smallest
$\sC (n,k,v)$, we have:

\begin{coro}\label{Cor2}
The size of the smallest banded asymmetric covering $\sD (n,1)$ is given by
$$C(n) = \sum_{i=0}^{[n/2]} C(n+1, n+1-2i, n-2i ) ~.$$
\end{coro}

Using the known values of $C(n,k,t)$ and the new values to be established in the next section (see Table \ref{T3}) we can determine $C(n)$ exactly for $n \le 11$.
These values are given in Table \ref{T1} and show that for $n \le 7$, banded asymmetric coverings are as good as any asymmetric coverings.

A more detailed investigation provides further information:

\begin{theorem}\label{thFE}
For lengths $n=1,2,3,5,6$ and $7$ an optimal asymmetric covering
is necessarily banded, and the corresponding covering
designs of length $1$ higher are unique.
At length $4$ there are four nonisomorphic minimal asymmetric covers, as shown in Table \ref{T3}, two banded and two non-banded.
\end{theorem}

\paragraph{Proof.}
By direct enumeration.
The details are omitted.~~~$\bsq$
\begin{table}[htb]
\caption{The four nonisomorphic minimal asymmetric coverings $\sD (4,1)$, of size $D(4,1) =6$.
(a) and (b) are banded, (c) and (d) are not.}
$$
\begin{array}{cccclcccclcccclcccc}
1 & 1 & 1 & 1 & ~~~~~~ &
1 & 1 & 1 & 1 & ~~~~~~ &
1 & 1 & 1 & 1 & ~~~~~~ &
1 & 1 & 1 & 1 \\
1 & 1 & 1 & 0 &&
1 & 1 & 1 & 0 &&
1 & 1 & 1 & 0 &&
1 & 1 & 1 & 0 \\
1 & 0 & 0 & 1 &&
1 & 1 & 0 & 1 &&
1 & 1 & 0 & 1 &&
1 & 1 & 0 & 1 \\
0 & 1 & 0 & 1 &&
0 & 0 & 1 & 1 &&
0 & 0 & 1 & 1 &&
0 & 0 & 1 & 1 \\
0 & 0 & 1 & 1 &&
1 & 1 & 0 & 0 &&
1 & 0 & 0 & 0 &&
0 & 1 & 0 & 1 \\
0 & 0 & 0 & 0 & &
0 & 0 & 0 & 1 & &
0 & 1 & 0 & 0 &&
1 & 0 & 0 & 0 \\
\multicolumn{4}{c}{{\rm (a)}}
&& \multicolumn{4}{c}{{\rm (b)}}
&& \multicolumn{4}{c}{{\rm (c)}}
&& \multicolumn{4}{c}{{\rm (d)}}
\end{array}
$$
\label{T3}
\end{table}

At length 7 the unique optimal (and banded) asymmetric covering can be found by deleting any coordinate from the following set of 31 vectors of length
8: $1^8$;
$\{0^2 1^2\}1^4$ (6);
$1^4 \{0^2 1^2 \} (6)$;
the 14 vectors of the Steiner system
$S(3,4,8)$;
$1^2 0^6$, $0^2 1^2 0^4$, $0^4 1^2 0^2$, $0^6 1^2$.
\paragraph{Remark.}
The continuous linear programming problem for the banded case is easily solved, and has size exactly $2^{n+1} / (n+2)$, which is asymptotically
\beql{EqE100}
2^{n+1} \left( \frac{1}{n} - \frac{2}{n^2} + O \left( \frac{1}{n^3} \right) \right) ~,
\eeq
just slightly worse than (\ref{EqE99}).
\section{New values for covering numbers}
Let $N(n,k,k-1,M)$ denote the number of nonisomorphic minimal covering designs $\sC (n,k,k-1)$ of size $M$, where of course $M \ge C (n,k,k-1)$.
The main results of this section are shown in Tables \ref{TC1} and \ref{TC2}.
\begin{table}[htb]
\caption{Covering numbers $C(n,k,k-1)$ for $n \le 12$.  Starred entries are new; see text for missing entries.}
$$
\begin{array}{c|cccccccccc}
n/k & 2 & 3 & 4 & 5 & 6 & 7 & 8 & 9 & 10 & 11 \\ \hline
2 & 1 \\
3 & 2 & 1 \\
4 & 2 & 3 & 1 \\
5 & 3 & 4 & 4 & 1 \\
6 & 3 & 6 & 6 & 5 & 1 \\
7 & 4 & 7 & 12 & 9 & 6 & 1 \\
8 & 4 & 11 & 14 & 20 & 12 & 7 & 1 \\
9 & 5 & 12 & 25 & 30 & 30 & 16 & 8 & 1 \\
10 & 5 & 17 & 30 & 51* & 50 & 45 & 20 & 9 & 1 \\
11 & 6 & 19 & 47 & 66 & \alpha & 84* & 63 & 25 & 10 & 1 \\
12 & 6 & 24 & 57 & 113 & 132 & \beta & 126* & 84 & 30 & 11 \\
13 & 7 & 26 & 78 & \gamma & 245 & \delta & \epsilon & 185* & 112 & 36
\end{array}
$$
\label{TC1}
\end{table}

Table \ref{TC1} gives the values of $C(n,k,k-1)$ for $n \le 12$.
Starred entries are new, and we have also shown that
$$C(14,10,9) = 259 ~.$$
In every case the coverings achieving these bounds were already known, see \cite{LJCR} for references.
Our contribution has been to show that no smaller covering exist.
The five remaining gaps in Table \ref{TC1} are at $C(11,6,5)$
(where $96 \le \alpha \le 100$),
$C(12,7,6)$ ($165 \le \beta \le 176$),
$C(13,5,4)$ ($149 \le \gamma \le 157$),
$C(13,7,6)$ ($257 \le \delta \le 264$),
$C(13,8,7)$ ($269 \le \epsilon \le 297$),
the lower bounds being new, except for $\gamma$.

Table \ref{TC2} gives values of $N(n,k,k-1,M)$ together with a brief indication of how they were found.
An entry such as
$$11,8,7: ~ 63 (40);
~64(1193) (\mbox{based on $[10,7,6,46]$})
$$
indicates that there are 40 nonisomorphic
minimal coverings $\sC (11,8,7)$ of size 63,
1193 of size 64, and that the latter enumeration was based on examining all possible ways to extend minimal covering designs
$\sC (10,7,6)$ of size $\le 46$.
(For any $\sC (11,8,7)$ of size 63
must puncture to a covering $\sC (10,7,6)$ which contains a minimal $\sC (10,7,6)$ of size $\le 46$.)
The symbol $\bsq$ at the end of a line in the tables indicates that the enumeration of minimal covering designs $\sC (n,k,k-1)$ for these values of $n$ and $k$ is complete.

It is worth drawing attention to the gaps that occur
just above the parameters corresponding to the 
Steiner coverings $\sC (8,4,3)$, 
$\sC(10,4,3)$, $\sC(11,5,4)$, $\sC(12,6,5)$.
For example a $\sC(12,6,5)$ that does not contain the 132-block Witt
design must contain at least 137 elements.

We know of no earlier table of this type, although isolated values have been published.
For example de Caen et~al. \cite{CKRM} showed that
$N(9,5,4,30) =3$ and $N(10,6,5,50) =1$.
The Steiner triple systems $S(2,3,n)$ have been enumerated for $n \le 19$ \cite{Kas}:
this gives the number of optimal $\sC (n,3,2)$'s for $n \equiv 1$ or $3~(\bmod~6)$.
Also Steiner quadruple systems $S(3,4,n)$ have been enumerated for $n \le 15$
(see the survey article \cite{HP});
this gives the number of optimal $\sC (n,4,3)$ for $n \equiv 2$ or $4~(\bmod ~6)$.

To test isomorphism we generally used the isomorphism
subroutines in the Magma computer
algebra system \cite{Mag1}, \cite{Mag2}, \cite{Mag3}.

To compute the entries $N(n,k,k-1,M)$ in Table \ref{TC2} and to establish the new lower bounds implicit in Table \ref{TC1} we made use of two different branch-and-bound procedures.

The first procedure branched by selecting one of the 
uncovered $(k-1)$-subsets which had the fewest remaining $k$-subsets which could cover it.
Let $\{x_1, x_2, \ldots, x_l \}$ denote the $k$-subsets which could cover it.
The procedure
recursively considered the $l$ alternatives $\{\{x_i$ in the covering, $x_j$ not in the covering, for $1 \le j \le i-1 \}$, $1 \le i \le l \}$.
This branching continued until either every $(k-1)$-subset was covered,
or until the lower bound
\beql{Seq1}
\left\lceil \sum_{v' \in T} \frac{1}{\max_{v:v' \subset v} (|\{v''\in T: v'' \subset v \}|)} \right\rceil
\eeq
on the number of available $k$-subsets needed to cover the set $T$ of uncovered $(k-1)$-subsets showed that no covering of size $M$ could be obtained from the current branch.

The second procedure used solutions of the continuous LP problem (see Section 3) to guide it.
It branched on the $k$-subset $x$ whose corresponding variable was closest to $\frac{1}{2}$,
considering the alternatives ``$x$ in the covering'' 
and ``$x$ not in the covering''.
This branching continued until one of the following obtained:

(a)~every $k$-subset had been placed in or excluded from the covering,

(b)~some $(k-1)$-subset could no longer be covered, or

(c)~applying the following lemma, where $U$ is the set of available $k$-subsets
and $T$ is the set of uncovered $(k-1)$-subsets, showed that
no covering of size $M$ could be obtained from the
current branch.  The optimal choice of $w$ in this bound is given by the
solution to the dual of the linear programming relaxation; the program used
a discrete, exact approximation to these dual variables for its bound.

\begin{lemma}
Let $T$ and $U$ be finite sets equipped with a relation $x\subset y$
for $x\in T$, $y\in U$, and let
$$
w:T\to \RR
$$
be a function satisfying
\beql{wvalid}
\sum_{x \in T : x\subset y} w(x)\le 1 ~ \mbox{for all}~ y \in U ~.
\eeq
Then any $\sC \subseteq U$ covering $T$ satisfies the lower bound
\[
|\sC|\ge \left\lceil \sum_{x\in T} w(x) \right\rceil ~.
\]
\end{lemma}

\paragraph{Proof.}
Let $\sC \subseteq U$ cover $T$.  Then
\begin{eqnarray*}
|\sC| & = & \sum_{y \in \sC} 1\\
& \geq & \sum_{y \in \sC} \sum_{x \in T : x \subset y} w(x) \\
& \geq & \sum_{x \in T} w(x)
\end{eqnarray*}
where the first inequality is from (\ref{wvalid}), and the second is
because $\sC$ covers $T$.  Since $|\sC|$ is an integer, the result follows.~~~$\bsq$

\paragraph{Remark.}
The bound \eqn{Seq1} is the special case of the Lemma
in which $U$ is the set of all
available $k$-subsets, $T$ is the set of uncovered $(k-1)$-subsets, and
$$
w(v)=\frac{1}{\max_{v:v' \subset v} (|\{v''\in T: v'' \subset v \}|)} ~.
$$

%This bound is valid for any choice of $w_{v'}$ satisfying \eqn{Seq2},
%and the bound \eqn{Seq1} comes from a specific choice.
%The dual of the linear programming relaxation gives the $w_{v'}$ which maximize this bound, and the program used a discrete, exact approximation to these dual variables for its bound.

Because the second program uses a stronger bound, it searched smaller branch-and-bound trees, but since it solved the linear programming relaxation at each node, it took more time per node.
As a result, the first program was more efficient for ``easy'' problems, and the second program for ``difficult'' problems (roughly, those
in Table 4 in the region bounded by $n \ge 9$ and $4 \le k \le n-4$).

\begin{table}[htb]
\caption{Values of $N(n,k,k-1, M)$, the number of nonisomorphic minimal covering designs $\sC (n,k,k-1)$ of size $M$.
See text for further details.}
$$\begin{array}{ll}
4,2,1: & 2 (1), 3 (1) \quad\bsq \\
5,3,2: & 4 (1), 5 (1), 6 (1) \quad\bsq \\
6,4,3: & 6 (1), 7 (1), 8 (1), 10 (1) \quad\bsq \\
7,5,4: & 9 (1), 11 (2), 12 (2), 15 (1) \quad\bsq \\
8,6,5: & 12  (1), 13 (1), 15 (2), 16 (3), 17 (2), 21 (1) \quad\bsq \\
9,7,6: & 16  (1), 18 (1), 19 (2), 20 (3), 21 (2), 22 (3), 23 (3), 28 (1) \quad\bsq \\
10,8,7: & 20  (1), 21 (1), 24 (4), 25 (2), 26 (7), 27 (5), 28 (4), 29 (2), \\
&  30 (4), 36  (1) \quad\bsq \\
11,9,8: &25  (1), 27 (1), 29 (2) \mbox{(based on 10,8,7,21)} \\
12,10,9: & 30  (1), 31 (1), 34 (1), 35 (4) \mbox{(based on 11,9,8,29)} \\
13,11,10: & 36  (1), 38 (1), 40 (0) \mbox{(based on 12,10,9,31)} \\
~ \\
5,2,1: & 3  (1), 4 (1) \quad\bsq \\
6,3,2: & 6  (1), 7 (5), 8 (2), 10 (1) \quad\bsq \\
7,4,3: & 12  (4), 13 (57), 14 (139), 15 (24), 16 (6), 17 (1), 20 (1) \quad\bsq \\
8,5,4: & 20  (6), 21 (263), 22 (7340) \mbox{(based on 7,4,3,13)} \\
9,6,5: & 30  (2), 31 (16), 32 (863) \mbox{(based on 8,5,4,21)} \\
10,7,6: &   45  (20), 46 (609) \mbox{(based on 9,6,5,32)} \\
11,8,7: &   63  (40), 64 (1193) \mbox{(based on 10,7,6,46)} \\
12,9,8: &   84  (4), 85 (46), 86 (1423) \mbox{(based on 11,8,7,64)} 
\end{array}
$$
\label{TC2}
\end{table}

\addtocounter{table}{-1}
\begin{table}[htb]
\caption{(continued)}
$$\begin{array}{ll}
~ \\
6,2,1:  &    3  (1), 4 (2), 5 (1) \quad\bsq \\
7,3,2;  &    7 (1), 9 (14), 10 (40), 11 (60), 12 (7), 13 (1), 15 (1) \quad\bsq \\
8,4,3:  &   14 (1), 17 (13) \\
9,5,4:  &   30 (3), 31 (18), 32 (459) \mbox{(based on 8,4,3,17)} \\
10,6,5: &   50 (1), 52 (4), 53 (56), 54 (880) \mbox{(based on 9,5,4,32)} \\
11,7,6: &   84 (3), 85 (0) \mbox{(based on 10,6,5,54)} \\
12,8,7: &  126 (3), 127 (2), 128 (0) \mbox{(based on 11,7,6,84)} \\
13,9,8: &  185 (1), 186 (0) \mbox{(based on 12,8,7,127)} \\
14,10,9:&  259 (1) \mbox{(based on 13,9,8,185}) \\
~ \\
7,2,1:  &    4 (1), 5 (2), 6 (1) \quad\bsq \\
8,3,2:  &   11 (5), 12 (145) \mbox{(based on 7,2,1,4)} \\
9,4,3:  &   25 (77), 26 (5562), 27 (538969) \mbox{(based on 8,3,2,12)} \\
10,5,4: &   51 (40), 52 (3354) \mbox{(based on 9,4,3,26)} \\
~ \\
8,2,1:  &    4 (1), 5 (2), 6 (3), 7 (1) \quad\bsq \\
9,3,2:  &   12 (1), 13 (1), 14 (64) \mbox{(based on 8,2,1,4)} \\
10,4,3: &   30 (1), 33 (43) \mbox{(based on 9,3,2,13)} \\
11,5,4: &   66 (1), 70 (78) \mbox{(based on 10,4,3,30)} \\
12,6,5: &   132 (1), 137 (87) \mbox{(based on 11,5,4,66)} \\
~ \\
9,2,1:  &    5 (1), 6 (3), 7 (3), 8 (1) \quad\bsq \\
10,3,2: &   17 (58) \mbox{(based on 9,2,1,5)} \\
11,4,3: &   47 (\ge 95970) \mbox{(based on 10,3,2,17)} \\
\end{array}
$$
\label{TC3}
\end{table}


\clearpage
\begin{thebibliography}{99}

\bibitem{Mag1}
W. Bosma and J. Cannon,
{\em Handbook of Magma Functions},
Sydney, May~22, 1995.

\bibitem{Mag2}
W. Bosma, J. Cannon and G. Mathews,
Programming with algebraic structures:
Design of the Magma language, in
{\em Proceedings of the 1994 International Symposium on Symbolic and Algebraic Computation},
M. Giesbrecht, Ed.,
Association for Computing Machinery, 1994, 52--57.

\bibitem{Mag3}
W. Bosma, J. Cannon and C. Playoust,
The Magma algebra system I:
The user language,
{\em J. Symb. Comp.}
{\bf 24} (1997), 235--265.


\bibitem{Andw}
A. E. Brouwer, J. B. Shearer,
N. J. A. Sloane and W. D. Smith,
A new table of constant weight codes,
{\em IEEE Trans. Inform. Theory},
{\bf 36} (1990), 1334--1380.

\bibitem{CKRM}
D. de Caen,
D. L. Kreher, S. P. Radziszowski and W. H. Mills,
On the covering of $t$-sets with $(t+1)$-sets:
$C(9,5,4)$ and $C(10,6,5)$,
{\em Discrete Math.}, {\bf 92} (1991), 65--77.

\bibitem{CKW}
D. de Caen, D. L. Kreher and J. Wiseman,
On constructive upper bounds for the Tur\'{a}n numbers
$T(n, 2r+1, 2r )$, pp. 277--280 of
{\em Nineteenth Southern Conference on Combinatorics, Graph Theory, and Computing $($Baton Rouge, 1988$)$}, Congress. Numer. {\bf 65} (1988).

\bibitem{CEK}
J. Cooper, R. Ellis and A. Kahng,
Directed binary covering codes, preprint, 2002.

\bibitem{CPLEX}
{\em CPLEX Manual},
CPLEX Organization Inc., Incline Village, Nevada, 1991.

\bibitem{FGK93}
R. Fourer, D. M. Gay and B. W. Kernighan,
{\em AMPL: A Modeling Language for Mathematical Programming},
Scientific Press, San Francisco, 1993.

\bibitem{Gir90}
G. R. Giraud,
Remarques sur deux probl\`{e}mes extr\'{e}maux,
{\em Discrete Math.},
{\bf 84} (1990), 319--321.

\bibitem{LJCR}
D. Gordon,
{\em La Jolla Covering Repository},
published electronically at
\htmladdnormallink{http://www.ccrwest.org/cover.html}{http://www.ccrwest.org/cover.html}.

\bibitem{HP}
A. Hartman and K. T. Phelps,
Steiner quadruple systems, pp. 205--240 of J. H. Dinitz and D. R. Stinson, eds.,
{\em Contemporary Design Theory}, Wiley, NY, 1992.

\bibitem{Kas}
P. Kaski and P. R. J. \"{O}sterg{\aa}rd,
The Steiner triple systems of order 19, preprint.


\bibitem{MS77}
F. J. MacWilliams and N. J. A. Sloane,
{\em The Theory of Error-Correcting Codes},
North-Holland, Amsterdam, 1977.

\bibitem{MM92}
W. H. Mills and R. C. Mullin,
Coverings and packings, pp. 371--399 of J. H. Dinitz and D. R. Stinson,
editors,
{\em Contemporary Design Theory:
A Collection of Surveys},
Wiley, NY, 1992.

\bibitem{NO99}
K. J. Nurmela and P. R. J. \"{O}sterg{\aa}rd,
\htmladdnormallink{New coverings of $t$-sets with $(t+1)$-sets}{http://www.tcs.hut.fi/Publications/papers/tp1.ps.Z},
{\em J. Combin. Des.},
{\bf 7} (1999), 217--226.

\bibitem{NO99a}
K. J. Nurmela and P. R. J. \"{O}sterg{\aa}rd,
New coverings of $t$-sets with $(t+1)$-sets: Appendix,
published electronically at
\htmladdnormallink{http://www.tcs.hut.fi/Publications/papers/table3.html}{http://www.tcs.hut.fi/Publications/papers/table3.html}.

\bibitem{OEIS}
N. J. A. Sloane,
{\em \htmladdnormallink{The On-Line Encyclopedia of Integer Sequences}{http://www.research.att.com/~njas/sequences}},
published electronically at www.research.att.com/$\sim$njas/sequences.

\bibitem{CRC}
D. R. Stinson, Coverings, pp. 260--265 of C. J. Colbourn and J. H. Dinitz,
editors,
{\em The CRC Handbook of Combinatorial Designs}, CRC Press,
Boca Raton, FL 1996.
\end{thebibliography}
\end{document}